\documentclass{article} 

\usepackage[usenames,dvipsnames]{pstricks} 
\usepackage{epsfig}
\usepackage{pst-grad} 
\usepackage{pst-plot} 
\usepackage[space]{grffile} 
\usepackage{etoolbox} 
\makeatletter 
\patchcmd\Gread@eps{\@inputcheck#1 }{\@inputcheck"#1"\relax}{}{} 
\makeatother

\usepackage{amsmath}
\usepackage{amssymb}
\usepackage{amsthm}
\usepackage[T1]{fontenc} 
\usepackage{graphicx}

\usepackage{marginnote}

\newcommand{\intav}[1]{\mathchoice {\mathop{\vrule width 6pt height 3 pt depth -2.5pt\kern -8pt \intop}\nolimits_{\kern -6pt#1}} {\mathop{\vrule width5pt height 3 pt depth -2.6pt \kern -6pt \intop}\nolimits_{#1}}{\mathop{\vrule width 5pt height 3 pt depth -2.6pt \kern -6pt\intop}\nolimits_{#1}} {\mathop{\vrule width 5pt height 3 pt depth-2.6pt \kern -6pt \intop}\nolimits_{#1}}}

\usepackage{ifthen}\usepackage{color}

\usepackage{setspace}

\newtheorem{thm}{Theorem}%
\newtheorem{cor}{Corollary}
\newtheorem{Definition}{Definition}
\newtheorem{Lemma}{Lemma}

\newtheorem{Proposition}{Proposition}
\newtheorem{Remark}{Remark}

\makeatletter\patchcmd{\env@cases}{1.2}{.8}{}{}\makeatother

\begin{document}

\title{Flipping regularity via the Harnack approach and applications to nonlinear elliptic problems}

\author{Diego R. Moreira and Edgard A. Pimentel}
	\date{\today} 
\maketitle

\begin{abstract}
\noindent We prove an abstract result ensuring that one-sided geometric control yields two-sided estimates for functions satisfying general conditions. Our findings resonate in the context of nonlinear elliptic problems, including supersolutions to fully nonlinear elliptic equations and functions in the De Giorgi class. Among the consequences of our abstract results are regularity estimates, and conditions for a continuous function to be in the class of viscosity solutions.  We also prove that one-sided geometric control yields $L^pL^\infty$-estimates. It provides a converse to the implication in the De Giorgi-Nash-Moser theory. 

\medskip

\medskip 

\end{abstract}

\vspace{.1in}

\section{Introduction}\label{sec_introduction}

Let $\Omega\subset\mathbb{R}^d$ be open and bounded. We consider functions $u:\Omega\to\mathbb{R}$ satisfying a modulus of continuity from below and prove an abstract result. Namely, under general conditions on the function, geometric control from below flips, yielding a modulus of continuity also from above. The abstract conditions ensuring this flipping are a weak Harnack inequality and an $L^pL^\infty$-estimate. The generality of such requirements encompasses a wide latitude of problems. For instance, elliptic equations in the divergence form, supersolutions to fully nonlinear equations, and minimizers and $Q$-minimizers in the calculus of variations. 

After establishing our abstract result, several consequences follow. In particular, we prove that functions satisfying a weak Harnack inequality and an $L^pL^\infty$-estimate are in the class of viscosity solutions. Also, supersolutions to fully nonlinear elliptic equations satisfying a geometric control from below (e.g., $\omega$-semiconvexity, $C^{0,\alpha}$ or $C^{1,\alpha}$-moduli of continuity) belong to suitable H\"older spaces depending on the integrability of the source term. 

Our results keep similarities with recent findings concerning the obstacle problem; see, for example, \cite{Andersson_Lindgren_Shahgholian_2015}. However, they are available in the fully nonlinear setting, and cover general supersolutions, in line with the Caffarelli-Kohn-Nirenberg-Spruck a priori estimate \cite{Caffarelli_Kohn_Nirenberg_Spruck_1985}.

The analysis of conditions under which a geometric control from below yields a modulus of continuity from above first appeared in the work of Luis Caffarelli \cite{Caffarelli_1999}; see also the series of lectures \cite{Caffarelli_2009}. In \cite{Caffarelli_1999}, the author examines the role of the Harnack inequality in the flipping of a given modulus of continuity, controlling the function from below. In that context, the author considers functions whose difference with respect to affine functions (or constants) satisfy a Harnack inequality. Under these conditions, and additional usual assumptions, he proves the function solves a fully nonlinear elliptic equation. An important application of the findings in \cite{Caffarelli_2009} is in the realm of homogenization; see, for instance \cite{Caffarelli_1999a}. 

One-sided control has also played a role in recent developments of the general theory of viscosity solutions. For instance, in \cite{Braga_Moreira_2018} the authors produce new (sharp) regularity estimates for semiconvex viscosity supersolutions. Their findings resonate in various contexts, including a connection with the celebrated Caffarelli-Kohn-Nirenberg-Spruck a priori estimate \cite{Caffarelli_Kohn_Nirenberg_Spruck_1985}.

In \cite{Braga_Figalli_Moreira_2019}, one-sided control is also a critical ingredient. In that paper, the authors consider fully nonlinear equations with unbounded ingredients. They prove two (optimal) regularity results for the convex envelope of supersolutions. Compare Theorems 2.6 and 2.8 in \cite{Braga_Figalli_Moreira_2019}. The importance of one-sided control, in terms of a semiconvexity condition, is clear in the regularity estimate in \cite[Theorem 2.9]{Braga_Figalli_Moreira_2019}.

We consider functions satisfying two conditions. The first one amounts to a weak Harnack inequality for non-negative $u:B_1\to\mathbb{R}$, centred at $x_0\in B_1$. That is, for $\rho>0$ and $x_0\in B_1$ such that $B_\rho(x_0)\subset B_1$, we have
\begin{equation}\label{eq_whint}\tag*{[WH]}
	\left(\intav{B_{\rho/2}(x_0)}u^\varepsilon{\rm{d}}x\right)^\frac{1}{\varepsilon}\leq C_{{\rm WH}}\left(\inf_{B_{\rho}(x_0)}u+\rho^\vartheta\chi_{x_0}(\rho)\right),
\end{equation}
where $\varepsilon>0$, $\chi$ is a non-negative function, and $\vartheta>0$ is a fixed exponent. Our second condition is an $L^pL^\infty$-estimate centred at $x_0\in B_1$
\begin{equation}\label{eq_lplinfint}\tag*{[$L^pL^\infty$]}
	\left\|u\right\|_{L^\infty(B_{\rho/2}(x_0))}\leq C_{p,\infty}\left[\left(\intav{B_{\rho}(x_0)}|u|^p{\rm{d}}x\right)^\frac{1}{p}+\sigma(\rho)\right],
\end{equation}
holding for some $p>0$, and some modulus of continuity $\sigma(\cdot)$. As before, $\rho>0$ and $x_0\in B_1$ are such that $B_\rho(x_0)\subset B_1$.

Our first main contribution is an abstract result. Suppose a function $u:B_1\to\mathbb{R}$ has a modulus of continuity from below. Suppose further that $u-\ell$ satisfies {\rm [WH]} and [$L^pL^\infty$] for every affine function $\ell$. Then the geometric control from below flips, yielding a two-sided control for the function. We state this result in the form of a theorem.

\begin{thm}[Flipping geometry I]\label{thm_main1}
Let $u:B_1\to \mathbb{R}$. Let $x_0\in B_{1/2}$. Suppose that $u-\ell$ satisfies {\rm [WH]} and {\rm [$L^p L^\infty$]} for every affine function $\ell$. Suppose further there exists a modulus of continuity $\gamma(\cdot)$ such that
\[
	\inf_{B_\rho(x_0)} \left(u-\ell_{x_0}\right)\geq -\gamma(\rho),
\]
for every $\rho\in(0,1/2)$. Then
\[
	\sup_{B_{\rho/4}(x_0)} \left(u-\ell_{x_0}\right)\leq C\left(\sigma(\rho)+\rho^\vartheta\chi_{x_0}(\rho)+\gamma(\rho)\right),
\]
where $C=C(C_{{\rm WH}}, C_{p,\infty})$.
\end{thm}

Notice that, if $\vartheta\geq1$, $\gamma(\rho)$ and $\sigma(\rho)$ are of order $o(\rho)$, and $\chi$ is of order $o(1)$, as $\rho\to 0$, then $u$ is differentiable at $x_0$. A consequence of the arguments in the proof of Theorem \ref{thm_main1} concerns functions $u:B_1\to\mathbb{R}$ such that $u-u(x_0)$ satisfies \ref{eq_whint} and \ref{eq_lplinfint} for suitable points $x_0\in B_1$.

Also, one can take $x_0\in B_1$ in the statement of Theorem \ref{thm_main1}. In that case, $\rho\in(0,\rho^*(x_0))$, with $\rho^*(x_0)\in(0,1-\|x_0\|)$. This is relevant when applying the Theorem \ref{thm_main1} in the context of interior regularity theory. In that case, for any $\Omega'\Subset B_1$, one must consider $\rho^*:=\rho^*(x)$, uniformly in $x\in\Omega'$.

\begin{cor}[Flipping geometry II]\label{cor_main1}
Let $u:B_1\to \mathbb{R}$. Let $x_0\in B_{1/2}$. Suppose that $u-c$ satisfies {\rm [WH]} and {\rm [$L^p L^\infty$]} for every constant $c\in\mathbb{R}$. Suppose further there exists a modulus of continuity $\gamma(\cdot)$ such that
\[
	\inf_{x\in B_\rho(x_0)} \big(u-u(x_0)\big)\geq -\gamma(\rho),
\]
for every $\rho\in (0,1/2)$. Then
\[
	\sup_{B_{\rho/4}(x_0)} \left(u-u(x_0)\right)\leq C\left(\sigma(\rho)+\rho^\vartheta\chi_{x_0}(\rho)+\gamma(\rho)\right),
\]
where $C=C(C_{{\rm WH}},C_{p,\infty})$.
\end{cor}

The abstract findings in Theorem \ref{thm_main1} and Corollary \ref{thm_main1} have important consequences for nonlinear elliptic problems, which we detail in the paper. The first implication of Theorem \ref{thm_main1} concerns the inclusion of a given function in a class of viscosity solutions. We prove that if $u-\ell$ satisfies \ref{eq_whint} and \ref{eq_lplinfint}, then $u$ belongs to a class of viscosity solutions. This class is characterized by ellipticity constants depending only on the dimension and the data of the problem. This is the subject of Theorem \ref{thm_whleliimplyvisc}; see Section \ref{subsec_whleliimplyvisc}.

\medskip

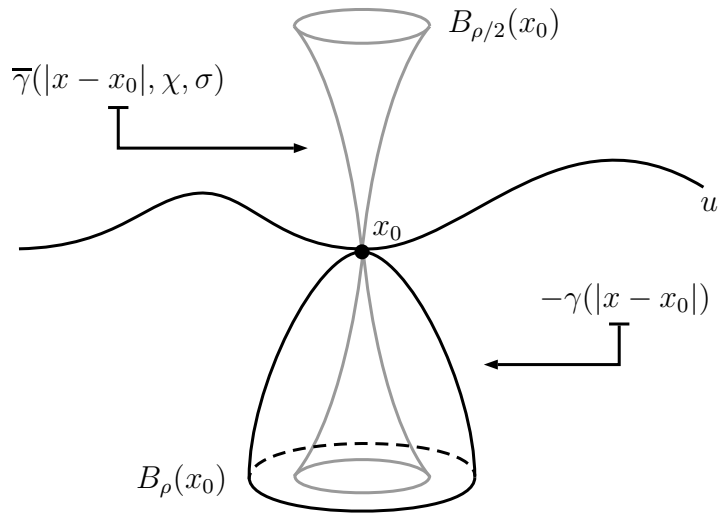
\begin{figure}[h!]
\center 
\psscalebox{1.0 1.0} 
{
\begin{pspicture}(0,-6.405)(9.41,0.345)
\definecolor{colour0}{rgb}{0.6,0.6,0.6}
\psbezier[linecolor=colour0, linewidth=0.04](4.656439,-2.9099588)(4.656439,-2.9099588)(4.8253813,-5.010877)(5.540633,-5.914655199199281)
\psbezier[linecolor=colour0, linewidth=0.04](4.656439,-2.9099588)(4.656439,-2.9099588)(4.4653864,-5.008982)(3.7406583,-5.905179344854892)
\psbezier[linecolor=colour0, linewidth=0.04](5.540633,-5.914655)(5.362215,-5.613712)(3.922365,-5.61413)(3.7407882,-5.913178289282679)
\psbezier[linecolor=colour0, linewidth=0.04](3.7485263,-5.9075665)(3.748848,-6.2027864)(5.540978,-6.208919)(5.5406337,-5.904139021071287)
\psbezier[linecolor=black, linewidth=0.04](0.08,-2.885)(1.4400709,-2.856904)(1.7657489,-2.1979268)(2.46,-2.145)(3.154251,-2.0920732)(3.5000708,-2.936904)(4.78,-2.885)(6.0599294,-2.833096)(7.334251,-0.9320732)(9.18,-2.065)
\psbezier[linecolor=black, linewidth=0.04](3.14,-5.925)(3.14,-4.725)(4.04,-2.925)(4.64,-2.925)(5.24,-2.925)(6.14,-4.725)(6.14,-5.925)
\psbezier[linecolor=black, linewidth=0.04](3.14,-5.925)(3.14,-6.525)(6.14,-6.525)(6.14,-5.925)
\psbezier[linecolor=black, linewidth=0.04, linestyle=dashed, dash=0.17638889cm 0.10583334cm](3.140146,-5.9459257)(3.14,-5.325)(6.14,-5.325)(6.139854,-5.904074283935638)
\psbezier[linecolor=colour0, linewidth=0.04](4.64,-2.925)(4.64,-2.925)(4.46,-0.825)(3.74,0.075)
\psbezier[linecolor=colour0, linewidth=0.04](4.64,-2.925)(4.64,-2.925)(4.82,-0.825)(5.54,0.075)
\psbezier[linecolor=colour0, linewidth=0.04](3.74,0.075)(3.92,-0.225)(5.36,-0.225)(5.54,0.075)
\psbezier[linecolor=colour0, linewidth=0.04](5.5399942,0.07978009)(5.54,0.375)(3.74,0.375)(3.7400062,0.0702199095488885)
\psdots[linecolor=black, dotsize=0.2](4.64,-2.925)
\rput[bl](9.14,-2.385){\large{$u$}}
\rput[bl](4.78,-2.765){\large{$x_0$}}
\rput[bl](7.0,-3.765){\large{$-\gamma(|x-x_0|)$}}
\rput[bl](0.0,-0.845){\large{$\overline{\gamma}(|x-x_0|,\chi,\sigma)$}}
\psline[linecolor=black, linewidth=0.04, tbarsize=0.07055555cm 5.0,arrowsize=0.05291667cm 2.0,arrowlength=1.4,arrowinset=0.0]{|*->}(1.4,-1.005)(1.4,-1.545)(3.92,-1.545)
\rput[bl](5.78,-0.145){\large{$B_{\rho/2}(x_0)$}}
\rput[bl](1.64,-6.185){\large{$B_\rho(x_0)$}}
\psline[linecolor=black, linewidth=0.04, tbarsize=0.07055555cm 5.0,arrowsize=0.05291667cm 2.0,arrowlength=1.4,arrowinset=0.0]{|*->}(8.06,-3.885)(8.06,-4.425)(6.26,-4.425)
\end{pspicture}
}
\caption{The modulus of continuity $\gamma$, touching $u$ at $x_0$ from below, flips to produce a modulus of continuity from above at $x_0$, denoted with $\overline{\gamma}$. However, the qualitative properties of $\overline{\gamma}$ are affected by two ingredients, namely: the factor $\sigma(\cdot)$ in the [$L^p L^\infty$] property, and $\chi$ in {\rm [WH]}. As a result, the one-side control encoded by $\gamma$ yields two-side information on the function. The latter may deteriorate, depending on the effects of $\sigma$ and $\chi$ appearing in the abstract conditions {\rm [WH]} and [$L^p L^\infty$].}\label{fig_changed}
\end{figure}

\medskip

After establishing this inclusion, we turn to the regularity theory. This important class of problems has known many recent advances, as several authors have worked in detail in various directions. The range of recent developments includes boundary regularity and differentiability \cite{Braga_Gomes_Moreira_Wang_2020,Ma_Moreira_Wang_2017,Silvestre_Sirakov_2014}, improved regularity under differentiability conditions \cite{Savin_2007,Armstrong_Silvestre_Smart_2012}, perturbative and approximation methods \cite{Teixeira_2014,Teixeira_Urbano_2014,Pimentel_2022}, non-uniformly elliptic problems \cite{Beck_Mingione_2020,Imbert_Silvestre_2016,DeFilippis_Mingione_2020}, regularity for elliptic variational problems \cite{Colombo_Mingione_2015,Colombo_Mingione_2015a,DeFilippis_Mingione_2020a,DeFilippis_Mingione_2021}, and estimates via potential methods \cite{Daskalopoulos_Kuusi_Mingione_2014,Kuusi_Mingione_2014a,Kuusi_Mingione_2014b,Kuusi_Mingione_2013,Kuusi_Mingione_2012,Mingione_2018} to name just a very few references.

 As concerns regularity theory, the Theorem \ref{thm_main1} and Corollary \ref{cor_main1} are consequential. Here, we work in the fully nonlinear setting, and the function $\chi_{x_0}$ in {\rm [WH]} becomes the norm of a right-hand side $f:B_1\to\mathbb{R}$ in suitable Lebesgue spaces. Let $u\in C(B_1)$, and $d/2\leq p_0<d$ be the exponent such that the Aleksandrov-Bakelman-Pucci estimate holds for $(\lambda,\Lambda)$-elliptic equations whose right-hand side is in $L^p(B_1)$, for $p_0<p$. If $u$ is a semiconvex viscosity supersolution in the presence of a right-hand side $f\in L^p(B_1)$, for $p_0<p<d$, Corollary \ref{cor_main1} implies $u\in C^{2-\frac{d}{p}}_{\rm loc}(B_1)$. The result also provides estimates and extends results in \cite{Braga_Figalli_Moreira_2019} to the range $p_0<p<d$. This fact is the content of the Theorem \ref{thm_hcer}. 

We observe that semiconvexity is a global property. Conversely, a uniform $C^{0,\alpha}$ or $C^{1,\alpha}$-modulus of continuity from below prescribes only a local geometric control. Our findings resonate also in this setting. We prove that viscosity supersolutions with a $C^{0,\alpha}$-modulus of continuity from below are H\"older continuous, if the right-hand side is in $L^p$, for $p_0<p$. Stemming from Corollary \ref{cor_main1}, this result also yields estimates.

When it comes to flipping a $C^{1,\alpha}$-modulus of continuity, we resort to Theorem \ref{thm_main1}. This is because to access information at the level of the gradient, we must subtract a tangent plane from the function under analysis. We consider supersolutions with a $C^{1,\alpha}$-modulus of continuity from below and right-hand side in $L^p$, for $p>d$. Here, we prove $C^{1,\beta}$-estimates, for some $\beta\in(0,1)$ universal. As before, our findings include estimates in terms of the usual quantities.

A further instance where Corollary \ref{cor_main1} is consequential concerns the so-called De Giorgi class; see, for instance, \cite{DiBenedetto_Trudinger_1984,DiBenedetto_Gianazza_2016}. Typically it comprises functions satisfying Caccioppoli inequalities, at any side and any truncation level. It includes solutions to elliptic equations in the divergence form as well as minimizers and $Q$-minimizers in the calculus of variations. Therefore the relevance of abstract results holding at the full generality of the De Giorgi class. In \cite{DiBenedetto_Gianazza_2016} the authors prove that functions in this class satisfy two properties; namely: a weak Harnack inequality and an $L^pL^\infty$-estimate. Hence, Corollary \ref{cor_main1} and its consequences are available for functions in the De Giorgi class.

We examine the causality involving the $L^pL^\infty$-estimates and H\"older continuity. Indeed, it is classical in the literature that the former operates as an ingredient in the proof of the latter. We establish a converse to this fact: a function satisfying a modulus of continuity from below satisfies an $L^pL^\infty$-estimate. Our first result in this direction concerns functions with a modulus of continuity of the zeroth order from below. That is, if a translation of the function satisfies a modulus of continuity from below, an $L^pL^\infty$-estimate is available. In addition, we consider functions satisfying a modulus of continuity from below at the level of the gradient. In this case, the assumption is that the function subtracted a tangent plane has some geometric control from below. Here, an $L^pL^\infty$-estimate is also available. The strategies used in both proofs differ, as in the second case the control of the gradient is pivotal.

The remainder of this manuscript is organized as follows. Section \ref{sec_prelim} gathers preliminary facts and basic material we use in the paper. The proofs of Theorem \ref{thm_main1} and Corollary \ref{cor_main1} are the subject of Section \ref{sec_flipping}. In Section \ref{sec_lplinflocal} we derive $L^pL^\infty$-estimates for functions satisfying one-sided geometric control. Finally, the paper details various consequences of the abstract results to nonlinear elliptic problems; this is the content of Section \ref{sec_cons}.

\section{Preliminaries}\label{sec_prelim}

In what follows, we gather preliminary material used throughout the paper. We start with a rigorous description of the conditions appearing in our abstract results.
\begin{enumerate}
\item[{\rm [WH]}] Let $u:B_1\to\mathbb{R}$ be nonnegative. Suppose $\chi$ is a non-negative function, and let $\vartheta>0$. We say $u$ satisfies {\rm [WH]} centred at $x_0\in B_1$ if there exists $\varepsilon>0$ such that
\[
	\left(\intav{B_{\rho/2}(x_0)}u^\varepsilon{\rm{d}}x\right)^\frac{1}{\varepsilon}\leq C_{{\rm WH}}\left(\inf_{B_{\rho/2}(x_0)}u+\rho^\vartheta\chi_{x_0}(\rho)\right),
\]
for every $\rho>0$ and $x_0\in B_1$ whenever $B_\rho(x_0)\subset B_1$.
\end{enumerate}

\begin{enumerate}
\item[{[$L^p L^\infty$]}] 
Let $u:B_1\to\mathbb{R}$ and $\sigma:[0,+\infty)\to[0,+\infty)$. We say $u:B_1\to\mathbb{R}$ satisfies an $L^p L^\infty$-estimate centred at $x_0\in B_1$ if there exists $C_{p,\infty}>0$ such that
\[
	\left\|u\right\|_{L^\infty(B_{\rho/2}(x_0))}\leq C_{p,\infty}\left[\left(\intav{B_{\rho}(x_0)}|u|^p{\rm{d}}x\right)^\frac{1}{p}+\sigma(\rho)\right]
\]
for some $p>0$, and every $\rho>0$ and $x_0\in B_1$ whenever $B_\rho(x_0)\subset B_1$.
\end{enumerate}

Our reasoning combines {\rm [WH]} and [$L^p L^\infty$]. A necessary condition for both inequalities to match sits in the exponents $\varepsilon>0$ and $p>0$, which must agree in {\rm [WH]} and [$L^p L^\infty$]. Although there is no a priori reason for this to be the case, the next lemma ensures this is always possible.

\begin{Lemma}[$L^pL^q$-interpolation]\label{lem_eptop}
Let $u:B_1\to\mathbb{R}$ and take $\rho>0$ and $x_0\in B_1$ such that $B_\rho(x_0)\subset B_1$. Suppose that $u$ satisfies
\[
	\left\|u\right\|_{L^\infty(B_{\rho/2}(x_0))}\leq C\left[\left(\intav{B_{\rho}(x_0)}|u|^p{\rm{d}}x\right)^\frac{1}{p}+\sigma(\rho)\right],
\]
for some $p>0$. Then there exists $\overline C>0$, depending only on $p$, $q$ and the dimension $d$, such that 
\[
	\left\|u\right\|_{L^\infty(B_{\rho/2}(x_0))}\leq \overline C\left[\left(\intav{B_{\rho}(x_0)}|u|^q{\rm{d}}x\right)^\frac{1}{q}+\sigma(\rho)\right]
\]
for every $q>0$.
\end{Lemma}
For the proof of Lemma \ref{lem_eptop}, we refer the reader to \cite[Remark 4.4]{Kinnunen_2001}. Part of our analysis concerns the properties of viscosity solutions to fully nonlinear equations. To properly state our findings, we gather a few notions and facts in this realm.

Let $0<\lambda\leq\Lambda$ be fixed constants and denote with $S(d)\sim\mathbb{R}^\frac{d(d+1)}{2}$ the space of symmetric matrices of order $d$. We define the Pucci extremal operators $\mathcal{M}^\pm_{\lambda,\Lambda}:S(d)\to\mathbb{R}$ as
\[
	\mathcal{M}^+_{\lambda,\Lambda}(M):=\Lambda\sum_{e_i>0}e_i+\lambda\sum_{e_i<0}e_i
\]
and
\[
	\mathcal{M}^-_{\lambda,\Lambda}(M):=\lambda\sum_{e_i>0}e_i+\Lambda\sum_{e_i<0}e_i,
\]
where $e_1,e_2,\ldots,e_d$ stand for the eigenvalues of $M$. We notice that any symmetric matrix $M$ can be written as $M=M^+-M^-$, where $M^+$ stands for the positive part of $M$, and $M^-$ stands for its negative part. For convenience, we recall a definition of the extremal operators in terms of those quantities; indeed,
\[
	\mathcal{M}^+_{\lambda,\Lambda}(M)=\Lambda{\rm Tr}(M^+)-\lambda{\rm Tr}(M^-)
\]
and
\[
	\mathcal{M}^-_{\lambda,\Lambda}(M)=\lambda{\rm Tr}(M^+)-\Lambda{\rm Tr}(M^-).
\]

Once these operators are available, one introduces the class of viscosity solutions.

Let $0<\lambda\leq\Lambda$ be fixed constants, and $f\in L^p(B_1)\cap C(B_1)$. The class of viscosity subsolutions $\underline{S}(\lambda,\Lambda,f)$ gathers the functions $u\in C(B_1)$ satisfying 
\begin{equation}\label{eq_cemut}
	\mathcal{M}^+_{\lambda,\Lambda}(D^2u)\geq f
\end{equation}
in the viscosity sense. Similarly, the class of viscosity supersolutions $\overline{S}(\lambda,\Lambda,f)$ comprises every $u\in C(B_1)$ such that 
\begin{equation}\label{eq_marco}
	\mathcal{M}^-_{\lambda,\Lambda}(D^2u)\leq f
\end{equation}
in the viscosity sense. The class of viscosity solutions $S(\lambda,\Lambda,f)$ intersects both, i.e.,
\[
	S(\lambda,\Lambda,f):=\underline{S}(\lambda,\Lambda,f)\cap\overline{S}(\lambda,\Lambda,f).
\]

Meanwhile, the inequalities in \eqref{eq_cemut} and \eqref{eq_marco} imply inequalities in terms of matrix norms. In fact, for $M\in S(d)$,
\[
	\mathcal{M}^-_{\lambda,\Lambda}(M)\leq C_0
\]
yields
\[
	\left\|M^+\right\|\leq\frac{d\Lambda}{\lambda}\left\|M^-\right\|+ \frac{C_0}{\lambda}.
\]
Similarly, if 
\[
	\mathcal{M}^+_{\lambda,\Lambda}(M)\geq C_0
\]
one obtains
\[
	\left\|M^-\right\|\leq\ \frac{d\Lambda}{\lambda}\left\|M^+\right\|-\frac{C_0}{\lambda}.
\]
This observation motivates an alternative definition of the viscosity class. It relies on inequalities relating to the positive and negative parts of $D^2P$, where $P$ is a paraboloid tested against the solutions.
\begin{Definition}[Class of viscosity solutions]\label{def_visclass}
Let $L_1,L_2>0$ be constants depending solely on $\lambda$, $\Lambda$, and the dimension $d$. We say that $u\in C(B_1)$ is in the class $\underline{\mathcal{S}}(L_1,L_2,C_0)$ if, whenever a paraboloid $P$ touches $u$ from above at $x_0\in B_1$, we have
\[
	\|(D^2P(x_0))^-\|\leq L_1\left\|(D^2P(x_0))^+\right\|+ L_2\cdot C_0.
\]
Similarly, we say that $u\in C(B_1)$ is in the class $\overline{\mathcal{S}}(L_1,L_2,C_0)$ if whenever a paraboloid $P$ touches $u$ from below at $x_o\in B_1$, we have
\[
	\|(D^2P(x_0))^+\|\leq L_1\left\|(D^2P(x_0))^-\right\|-L_2\cdot C_0.
\]

\end{Definition}

We proceed with a lemma comparing the opening of paraboloids that are ordered in some neighbourhood. Before stating the lemma, we introduce some notation. Let $M\in \mathbb{R}^{d^2}$be a square matrix. We denote with $\lambda_{{\rm min}}(M)$ the smallest eigenvalue of $M$; with $\lambda_{{\rm max}}(M)$ we denote the largest eigenvalue of $M$.

\begin{Lemma}[Comparable openings]\label{lem_curvature}
Let $P(x):=A|x-x_0|^2$, where $A\geq 0$. Let $Q$ be a paraboloid such that $P$ touches $Q$ from above at some $x_0\in\mathbb{R}^d$. Then
\[
	A\geq \|(D^2Q(x_0))^+\|.
\]
\end{Lemma}
\begin{proof}
We start by supposing $(D^2Q(x_0))^+\neq0$, as otherwise the statement of the lemma trivially holds. In this case,
\[
	\lambda_{\rm max}(D^2Q(x_0))=\|(D^2Q(x_0))^+\|.
\]
Without loss of generality, suppose $x_0\equiv 0$ and the affine part of $Q$ is identically zero. The latter assumption is not restrictive, as the argument focuses on the second derivative of $Q$. Furthermore, let $\mathcal{B}:=\{e_1,\ldots,e_d\}$ denote an orthonormal basis such that 
\[
	Q(x)=\sum_{\lambda_i>0}\lambda_ix_i^2+\sum_{\lambda_i<0}\lambda_ix_i^2,
\]
where $Q(e_i)=\lambda_ie_i$. Denote with $e_{\rm max}$ the eigenvector associated with the eigenvalue $\lambda_{\rm max}(D^2Q(x_0))$. Because $P$ touches $Q$ from above in $B_\delta$, we have
\[
	A\left|\frac{\delta}{2}e_{\rm max}\right|^2=P\left(\frac{\delta}{2}e_{\rm max}\right)\geq Q\left(\frac{\delta}{2}e_{\rm max}\right)=\lambda_{\rm max}\left(\frac{\delta}{2}\right)^2.
\]
Hence
\[
	A\geq \|(D^2Q(x_0))^+\|
\]
and the proof is complete.
\end{proof}

We proceed with the definition of $\omega$-semiconvexity.

\begin{Definition}[$\omega$-semiconvex function]\label{def_omegasc}
Let $\omega:[0,+\infty)\to[0,+\infty)$ be a modulus of continuity. We say that $u:B_1\to\mathbb{R}$ is $\omega$-semiconvex if
\[
	u(tx+(1-t)y)\leq tu(x)+(1-t)u(y)+t(1-t)|x-y|\omega(|x-y|),
\]
for every $x,y\in B_1$ and every $t\in[0,1]$.
\end{Definition}

\begin{Remark}[Semiconvexity and affine translations]\label{rem_semiconvex}
Let $u:B_1\to\mathbb{R}$. Let $\ell(x):=\ell_0+\ell_1\cdot x$, where $\ell_0\in\mathbb{R}$ and $\ell_1\in\mathbb{R}^d$. Set $v:=u+\ell$. If $u$ is $\omega$-semiconvex, so is $v$. In fact,
\begin{align*}
	v(tx+(1-t)y)&=u(tx+(1-t)y)+t\ell_0+t\ell_1\cdot x+(1-t)\ell_0+(1-t)\ell_1\cdot y\\
		&\leq t(u(x)+\ell(x))+(1-t)(u(y)+\ell(y))\\
			&\quad+t(1-t)|x-y|\omega(|x-y|)\\
		&=tv(x)+(1-t)v(y)+t(1-t)|x-y|\omega(|x-y|),
\end{align*}
for every $x,y\in B_1$ and every $t\in[0,1]$.
\end{Remark}

We proceed with a lemma ensuring that $\omega$-semiconvex functions satisfy [$L^p L^\infty$], with a precise contribution from the modulus of continuity $\omega$.

\begin{Lemma}[$L^p L^\infty$-estimate for $\omega$-semiconvex functions]\label{lem_BFMprop82b}
Let $u\in L^1(B_1)$ be an $\omega$-semiconvex function. There exists $C>0$ such that 
\[
	\sup_{B_{\rho/2}}|u|\leq C\left(\left(\intav{B_\rho}|u|^p{\rm d}x\right)^\frac{1}{p}+\rho\omega(\rho)\right)
\]
for every $\rho\in(0,1]$ and every $p\geq 1$. In particular, $C=C(p,d,\omega(1))$.
\end{Lemma}

It follows from Lemma \ref{lem_BFMprop82b} that $\omega$-semiconvex functions satisfy [$L^p L^\infty$] with $\sigma(t):=t\omega(t)$. For a proof of this result, we refer the reader to \cite[Proposition 8.2, item (b)]{Braga_Figalli_Moreira_2019}. In the sequel, we include information on the subdifferentials of $\omega$-semiconvex functions. Given an $\omega$-semiconvex function $u:B_1\to \mathbb{R}$, its subdifferential at $x\in B_1$ is denoted with $\partial_\omega u(x)$, and defined as
\[
	\partial_\omega u(x):=\left\lbrace P\in\mathbb{R}^d\;|\;u(y)\geq u(x)+\left\langle P,y-x\right\rangle-|y-x|\omega(|y-x|),\;\forall\,y\in B_1\right\rbrace.
\]

\begin{Lemma}[Subdifferential estimates]\label{lem_BMProp51}
Let $u:B_1\to\mathbb{R}$ be a bounded function. Suppose $u$ is $\omega$-semiconvex. Then $\partial_\omega u(x)$ is non-empty, compact, and convex, for every $x\in B_1$. In addition, let $K\subset B_1$ be a compact set. Hence,
\[
	\sup_{x\in K}\sup_{P\in\partial_\omega u(x)}|P|\leq 2\left(\frac{\|u\|_{L^\infty(B_1)}}{{\rm dist}(K,\partial B_1)}+\omega\big({\rm dist}(K,\partial B_1)\big)\right).
\]
\end{Lemma}
We notice this result holds for any convex domain $\Omega\subset\mathbb{R}^d$ instead of $B_1$. For a proof of Lemma \ref{lem_BMProp51}, we refer the reader to \cite[Proposition 5.1]{Braga_Moreira_2018}. We continue with an inequality for bounded nonnegative functions defined on intervals.

\begin{Lemma}\label{lem_giusti}
Fix $\rho>0$ and let $\Phi:[\rho,R]\to\mathbb{R}$ be a bounded nonnegative function. Suppose that for $\rho\leq t<s\leq R$ we have
\[
	\Phi(t)\leq\theta\Phi(s)+\frac{C_1}{(s-t)^\alpha}+C_2,
\]
where $\theta\in(0,1)$, $C>0$, and $\alpha>0$. Then
\[
	\Phi(\rho)\leq C(\alpha,\theta)\left(\frac{C_1}{(R-\rho)^\alpha}+C_2\right).
\]
\end{Lemma}
For a proof of Lemma \ref{lem_giusti}, we refer the reader to \cite[Lemma 6.1]{Giusti_2003}. 

\section{Flipping geometry: an abstract result}\label{sec_flipping}

In the sequel we combine {\rm [WH]} and [$L^p L^\infty$] to establish our abstract results; we start with the proof of the Theorem \ref{thm_main1}.

\begin{proof}[Proof of Theorem \ref{thm_main1}]
We argue by combining {\rm [WH]} and {\rm [$L^p L^\infty$]}, and resorting to previous results. Start by fixing $x_0\in B_1$ and $\rho>0$ such that $B_\rho(x_0)\subset B_1$; define the auxiliary function $v:B_\rho\to\mathbb{R}$ as
\[
	v(x):=u(x)-\ell_{x_0}(x)+\gamma(\rho),
\]
for some affine function $\ell_{x_0}$ such that $\ell_{x_0}(x_0)=u(x_0)$. Clearly, $v$ is nonnegative. By assumption, $u-\ell_{x_0}$ satisfies {\rm [WH]} centred at $x_0\in B_1$; hence,
\begin{align}\label{eq_thm1_1}
	\nonumber\left(\intav{B_{\rho/2}(x_0)}v^\varepsilon\right)^\frac{1}{\varepsilon}&\leq C_{{\rm WH}}\left(\inf_{B_{\rho/2}(x_0)}v+\rho^\vartheta\chi_{x_0}(\rho)\right)\\
		\nonumber&\leq C_{{\rm WH}}\left(v(x_0)+\rho^\vartheta\chi_{x_0}(\rho)\right)\\
				&= C_{{\rm WH}}\left(\gamma(\rho)+\rho^\vartheta\chi_{x_0}(\rho)\right).
\end{align}
Moreover, the fact that $u$ and $v$ differ by an affine function builds upon {\rm [$L^p L^\infty$]} to yield
\begin{equation*}\label{eq_thm1_2}
	\left\|v\right\|_{L^\infty(B_{\rho/4}(x_0))}\leq C_{p,\infty}\left[\left(\intav{B_{\rho/2}(x_0)}v^p{\rm{d}}x\right)^\frac{1}{p}+\sigma(\rho)\right].
\end{equation*}
Because of Lemma \ref{lem_eptop}, the former inequality becomes
\begin{equation}\label{eq_thm1_3}
	\left\|v\right\|_{L^\infty(B_{\rho/4}(x_0))}\leq C\left[\left(\intav{B_{\rho/2}(x_0)}v^\varepsilon{\rm{d}}x\right)^\frac{1}{\varepsilon}+\sigma(\rho)\right],
\end{equation}
for the same $\varepsilon>0$ as in {\rm [WH]}. Hence, by combining \eqref{eq_thm1_1} and \eqref{eq_thm1_3} we obtain
\begin{equation}\label{eq_thm1_4}
	\sup_{B_{\rho/4}(x_0)}v\leq C\left[C_{{\rm WH}}\left(\gamma(\rho)+\rho^\vartheta\chi_{x_0}(\rho)\right)+\sigma(\rho)\right].
\end{equation}
On the other hand, the definition of $v$ implies
\begin{equation}\label{eq_thm1_5}
	\sup_{B_{\rho/4}(x_0)}u(x)-\ell_{x_0}(x)\leq \sup_{B_{\rho/4}(x_0)}v.
\end{equation}
Finally, combine \eqref{eq_thm1_4} and \eqref{eq_thm1_5}, and notice that $0<\rho\ll1$ was taken arbitrarily to complete the proof.
\end{proof}

\begin{Remark}\label{rem_rhopsycho}\normalfont
In some applications, the $L^pL^\infty$-estimate may appear prescribed in balls whose radii differ from the ones in {\rm [$L^p L^\infty$]}, as in Propositions \ref{prop_lplinftycalpha} and \ref{prop_lplinftyc1alpha}. In this case, Theorem \ref{thm_main1} and Corollary \ref{cor_main1} still hold true, under minor adjustments in the passage from \eqref{eq_thm1_2} to \eqref{eq_thm1_3}.
\end{Remark}

As a consequence of the Theorem \ref{thm_main1}, we infer the flipping of moduli of continuity at the level of the function. That is, in the case where $u-u(x_0)$ satisfies {\rm  [WH]} and {\rm [$L^p L^\infty$]} for $x_0\in B_1$. This is the content of Corollary \ref{cor_main1}, whose proof we discuss in the sequel.

\begin{proof}[Proof of Corollary \ref{cor_main1}]
The proof follows along the same lines as in the proof of Theorem \ref{thm_main1}. It suffices to consider the constant affine function $\ell\equiv u(x_0)$.
\end{proof}

\begin{Remark}[Gradient regularity vis-a-vis modulus of continuity for the function]\label{rmk_affinevavconstant}\normalfont
We notice a fundamental difference between Theorem \ref{thm_main1} and Corollary \ref{cor_main1}, regarding the regularity of the function $u$. Indeed, for a truly affine map $\ell(x):={\rm a}+{\rm b}\cdot x$, with ${\rm b}\neq 0$, we flip information for the function subtracted a tangent plane. As a consequence, the analysis yields regularity \emph{at the level of the gradient}. Conversely, if we are in the context of Corollary \ref{cor_main1}, where $\ell(x):=a$, the flipping yields a modulus of continuity \emph{for the function}.  
\end{Remark}

\begin{Remark}[One-sided regularity implies two-sided regularity]\label{rmk_twosided}\normalfont We have chosen to present Theorem \ref{thm_main1} and Corollary \ref{cor_main1} as flipping \emph{from below}. That is, given a modulus of continuity from below, the weak Harnack inequality and the local maximum principle produce a related modulus of continuity from above. However, by supposing that $\ell-u$ also satisfies {\rm [WH]} and [$L^p L^\infty$], the reverse implication holds. Given a modulus of continuity from above, the weak Harnack inequality combines with the $L^p L^\infty$-estimates and generates a modulus of continuity from below. Indeed, suppose $\ell-u$ satisfies [WH] and [$L^p L^\infty$]. Suppose further
\[
	\sup_{B_\rho(x_0)}\big(u(x)-\ell(x)\big)\leq \gamma(\rho)
\]
for some modulus of continuity $\gamma(\cdot)$, some affine function $\ell$ with $\ell(x_0)=u(x_0)$, and $0<\rho\ll1$ and $x_0\in B_1$ such that $B_\rho(x_0)\subset B_1$. Define $w:B_\rho(x_0)\to\mathbb{R}$ as
\[
	w(x):=\gamma(\rho)+\ell(x)-u(x)
\]
and notice that $w\geq 0$. Arguing as in the proof of Theorem \ref{thm_main1}, we obtain
\[
	\sup_{B_{\rho/4(x_0)}}\big(\ell(x)-u(x)\big)\leq C\left(\gamma(\rho)+\rho^\vartheta\chi_{x_0}(\rho)+\sigma(\rho)\right).
\]
This is the same as
\[
	\inf_{B_{\rho/4}(x_0)}\big(u(x)-\ell(x)\big)\geq -C\left(\gamma(\rho)+\rho^\vartheta\chi_{x_0}(\rho)+\sigma(\rho)\right).
\]
In conclusion, if both $\ell-u$ and $u-\ell$ satisfy {\rm [WH]} and {[$L^p L^\infty$]}, one-sided regularity implies two-sided regularity.
\end{Remark}

\begin{Remark}\normalfont
Let $L$ denote a purely second-order linear operator and suppose $g\in L^q(B_1)$. Let $u:B_1\to\mathbb{R}$ be such that $|Lu|\leq g$. Then $\ell-u$ and $u-\ell$ satisfy {\rm [WH]} and {[$L^\varepsilon L^\infty$]}. As a consequence, a one-sided modulus of continuity available for $u$ yields a two-sided regularity estimate. The same conclusion is available in the context of fully nonlinear operators, as long as $u$ satisfies
\[
	\mathcal{M}^-_{\lambda,\Lambda}(D^2u)\leq |g|\hspace{.2in}\mbox{and}\hspace{.2in}\mathcal{M}^+_{\lambda,\Lambda}(D^2u)\geq -|g|
\]
in $B_1$.
\end{Remark}

\section{Local $L^pL^\infty$-estimates}\label{sec_lplinflocal}

In this section, we prove that a modulus of continuity from below implies an $L^pL^\infty$-estimate. To some extent, our result reverses the usual implication available in the literature. Namely, that $L^pL^\infty$-estimates combined with energy inequalities yield H\"older-continuity. We mention, for instance, the role of $L^2L^\infty$-estimates in De Giorgi's solution of Hilbert's XIX problem \cite{DeGiorgi_1957}; see also \cite{DiBenedetto_Gianazza_2016, DiBenedetto_Trudinger_1984}.

More precisely, the results in the present section reverse this usual implication in the following sense. Fix $R>0$; by requiring $u\in L^\infty(B_R)$ to have a uniform modulus of continuity \emph{from below} in $B_{R/2}$, we prove an $L^pL^\infty$-estimate in $B_{R/16}$. Our findings include moduli of continuity from below both at the level of the function as well as at the level of first-order derivatives. We proceed with the former.

\begin{Proposition}[Local $L^pL^\infty$-estimates I]\label{prop_lplinftycalpha}
Fix $R>0$, and let $u\in L^\infty(B_R)$. Let $\sigma:[0,\infty)\to[0,\infty)$ be a modulus of continuity. Suppose
\[
	u(x)\geq u(y)-\sigma(|x-y|),
\]
for every $x\in B_{R/2}$, and every $y\in B_{R/8}(x)$. For every $p\in(0,\infty)$ there exists a positive constant $C=C(d,p)$ such that 
\[
	\|u\|_{L^\infty(B_{R/16})}\leq C\left(\left(\intav{B_R}|u|^p\right)^\frac{1}{p}+\sigma(R)\right).
\]
\end{Proposition}
\begin{proof}
We split the proof of the proposition into three steps.

\medskip

\noindent{\bf Step 1 - }We start by producing pointwise bounds for $u$ in $B_{R/2}$. Indeed, let $z\in B_{R/2}$; for any $y\in B_{R/8}(z)\subset B_R$ we have
\[
	u(y)\geq u(z)-\sigma(|y-z|)\geq u(z)-\sigma(R),
\]
since $\sigma(\cdot)$ is nondecreasing. As a consequence,
\begin{equation}\label{eq_firenze}
	u(z)\leq \intav{B_{R/8}(z)}u(y){\rm d}y+\sigma(R)\leq8^d\intav{B_{R}}|u(y)|{\rm d}y+\sigma(R).
\end{equation}

\medskip

\noindent{\bf Step 2 - }On the other hand, for $z,y\in B_{R/16}$, we have $|y-z|\leq R/8$. The monotonicity of $\sigma(\cdot)$ then implies
\begin{equation}\label{eq_mostly}
	u(y)\leq u(z)+\sigma(R).
\end{equation}

Let $\xi\in C^\infty_0(B_1)$ be such that $0\leq \xi\leq 1$, with
\[
	\xi(x)\equiv
		\begin{cases}
			1&\hspace{.3in}\mbox{in}\hspace{.15in}B_{R/32}\\
			0&\hspace{.3in}\mbox{in}\hspace{.15in}\mathbb{R}^d\setminus B_{R/16}.
		\end{cases}
\]
Multiplying both sides of \eqref{eq_mostly} by $\xi$ and integrating over $B_{R/32}$, we get
\[
	u(z)\int_{B_{R/16}}\xi(y){\rm d}y\geq-\left(\int_{B_{R/16}}|u(y)|{\rm d}y+\int_{B_{R/16}}\sigma(R){\rm d}y\right).
\]
Because 
\[
	\int_{B_{R/16}}\xi(y){\rm d}y\geq \frac{|B_{R/16}|}{2^d},
\]
the former inequality becomes
\begin{align}\label{eq_peixe}\nonumber
	u(z)&\geq-\frac{2^d}{|B_{R/16}|}\left(\int_{B_{R/16}}|u(y)|{\rm d}y+\int_{B_{R/16}}\sigma(R){\rm d}y\right)\\
		&\geq -32^d\left(\intav{B_R}|u(y)|{\rm d}y+\sigma(R)\right).
\end{align}
By combining \eqref{eq_peixe} and \eqref{eq_firenze}, we get
\[
	\|u\|_{L^\infty(B_{R/16})}\leq 32^d\left(\intav{B_R}|u(y)|{\rm d}y+\sigma(R)\right)
\]
and complete the proof in the case $p=1$. The monotonicity of the average yields the result  for every $p>1$ as well. It remains to verify the case $p\in (0,1)$.

\medskip

\noindent{\bf Step 3 - }If $p\in (0,1)$, the result follows from elementary inequalities combined with Lemma \ref{lem_giusti}. Let $0<R/16\leq t<s\leq R$ and take $x_0\in B_{R/16}$. Apply the statement of the proposition for $p=1$ to obtain
\[
	\begin{split}
		|u(x_0)|&\leq\sup_{B_{(s-t)/16}(x_0)}|u(x)|\\
			&\leq C\left(\frac{1}{(s-t)^d}\int_{B_s}|u(x)|{\rm d}x+\sigma(R)\right)\\
			&\leq C\left(\sup_{B_s}|u|^{1-p}\frac{1}{(s-t)^d}\int_{B_s}|u(x)|^p{\rm d}x+\sigma(R)\right).
	\end{split}	
\]
By applying Young's inequality with $\varepsilon\equiv1/2C$ and taking the supremum in $x_0\in B_t$, one gets
\[
	\sup_{B_t}|u(x)|\leq \frac{1}{2}\sup_{B_s}|u(x)|+C\left(\frac{2Cp(1-p)^\frac{1-p}{p}}{(s-t)^\frac{d}{p}}\left(\int_{B_s}|u(x)|^p{\rm d}x\right)^\frac{1}{p}+\sigma(R)\right).
\]
Now, Lemma \ref{lem_giusti} implies
\[
	\begin{split}
		\sup_{B_{R/16}}|u(x)|&\leq C(p,C)\left(\left(\frac{16}{15R}\right)^\frac{d}{p}\left(\int_{B_R}|u(x)|^p{\rm d}x\right)^\frac{1}{p}+\sigma(R)\right)\\
			&\leq C(p,C,d)\left(\left(\intav{B_R}|u(x)|^p{\rm d}x\right)^\frac{1}{p}+\sigma(R)\right)
	\end{split}
\]
and finishes the proof.
\end{proof}

Our next result derives an $L^pL^\infty$-estimate for functions satisfying a modulus of continuity at the level of first-order derivatives.

\begin{Proposition}[Local $L^pL^\infty$-estimates II]\label{prop_lplinftyc1alpha}
Let $u\in L^\infty(B_R)$ for some $R>0$ fixed, though arbitrary. Let $\sigma:[0,\infty)\to[0,\infty)$ be a modulus of continuity. Suppose that for every $y\in B_{R/2}$ there exists $P_y\in\mathbb{R}^d$ such that 
\begin{equation}\label{eq_belowc1a}
	u(x)\geq u(y)+\left\langle P_y,x-y\right\rangle-|x-y|\sigma(|x-y|),
\end{equation}
for every $x\in B_{R/8}(y)$. For every $p\in(0,\infty)$, there exists a positive constant $C=C(d,p)$ such that 
\[
	\|u\|_{L^\infty(B_{R/32})}\leq C\left(\left(\intav{B_R}|u|^p\right)^\frac{1}{p}+R\sigma(R)\right).
\]
\end{Proposition}
\begin{proof}
For ease of presentation, we split the proof into four steps. We start by supposing $u\in C^1(B_R)$ and resorting to \eqref{eq_belowc1a}.

\medskip

\noindent{\bf Step 1 - }Suppose $u\in C^1(B_R)$; hence $P_y\equiv Du(y)$. Using this fact, and averaging both sides of \eqref{eq_belowc1a} in $B_{R/8}(y)$, one obtains
\[
	\begin{split}
		\intav{B_{R/8}(y)}u(x){\rm d}x&\geq u(y)+\intav{B_{R/8}(y)}\left\langle Du(y),x-y\right\rangle{\rm d}x-R\sigma(R)\\
			&\geq u(y)+\left\langle Du(y),y-y\right\rangle-R\sigma(R)\\
			&=u(y)-R\sigma(R),
	\end{split}
\]
where the second inequality follows from the harmonicity of $\left\langle Du(y),x-y\right\rangle$ in $x$. Observe that
\[
	B_{R/8}(y)\subset B_{5R/8}\subset B_{3R/4}(y);
\]
hence,
\begin{equation}\label{eq_firstineqlinflp}
		u(y)\leq\intav{B_{R/8}(y)}u(x){\rm d}x+R\sigma(R)\leq 6^d\left(\intav{B_{3R/4}(y)}|u(x)|{\rm d}x+R\sigma(R)\right),
\end{equation}
for every $y\in B_{R/2}$. Next, we produce a lower bound for $u(y)$.

\medskip

\noindent{\bf Step 2 - }Let $x,y\in B_{R/16}\subset B_{R/2}$. Because $|x-y|\leq R/8$, we have $y\in B_{R/8}(x)$. Then
\begin{equation}\label{eq_gangemi}
	u(x)\geq u(y)+\left\langle Du(y),x-y\right\rangle-R\sigma(R),
\end{equation}
where we used the monotonicity of $\sigma$. As before, consider $\xi\in C^\infty_0(B_R)$ with $0\leq \xi\leq 1$, such that 
\[
	\xi(x)\equiv
		\begin{cases}
			1&\hspace{.3in}\mbox{in}\hspace{.15in}B_{R/32}\\
			0&\hspace{.3in}\mbox{in}\hspace{.15in}\mathbb{R}^d\setminus B_{R/16}.
		\end{cases}
\]
Suppose further
\[
	|D\xi|\leq \frac{C}{R},
\]
for some constant $C>0$. Multiply \eqref{eq_gangemi} by $\xi$ and integrate over $B_{R/16}$ to get
\begin{align}\label{eq_straciatella}\nonumber
	u(x)\int_{B_{R/16}}\xi(y){\rm d}y&\geq \int_{B_{R/16}}u(y)\xi(y){\rm d}y+\int_{B_{R/16}}\left\langle Du(y),x-y\right\rangle\xi(y){\rm d}y\\
			&\quad-R\sigma(R)\int_{B_{R/16}}\xi(y){\rm d}y.
\end{align}
We continue with the analysis of the term involving $Du$; because $\xi\equiv 0$ on $\partial B_{R/16}$, one obtains
\[
	\begin{split}
		0&=\int_{B_{R/16}}{\rm div}_y\big(u(y)\xi(y)(x-y)\big){\rm d}y\\
			&=\int_{B_{R/16}}\left\langle Du(y),x-y\right\rangle\xi(y){\rm d}y+\int_{B_{R/16}}u(y){\rm div}_y\big(\xi(y)(x-y)\big){\rm d}y.
	\end{split}
\]
Notice that 
\[
	\left|{\rm div}_y\big(\xi(y)(x-y)\big)\right|\leq C
\]
uniformly in $x,y\in B_{R/16}$, for some $C>0$ depending only on the dimension. Hence, the previous equality yields
\[
	\int_{B_{R/16}}\left\langle Du(y),x-y\right\rangle\xi(y){\rm d}y\leq C\int_{B_{R/16}}|u(y)|{\rm d}y. 
\]
Therefore, \eqref{eq_straciatella} becomes
\[
	\begin{split}
		u(x)\int_{B_{R/16}}\xi(y){\rm d}y&\geq -(1+C)\int_{B_{R/16}}u(y){\rm d}y-R\sigma(R)\int_{B_{R/16}}\xi(y){\rm d}y.	
	\end{split}
\]
Recalling that 
\[
	\int_{B_{R/16}}\xi(y){\rm d}y\geq \frac{|B_{R/16}|}{2^d},
\]
and arguing as before, one concludes that
\begin{equation}\label{eq_ferrero}
	u(x)\geq -24^d(C+1)\left(\intav{B_{R/16}}|u(y)|{\rm d}y+R\sigma(R)\right).
\end{equation}

\medskip

\noindent{\bf Step 3 - }By combining \eqref{eq_firstineqlinflp} and \eqref{eq_ferrero} we get
\[
	\left\|u\right\|_{L^\infty(B_{R/16})}\leq C\left(\intav{B_{R/16}}|u(y)|{\rm d}y+R\sigma(R)\right)
\]
and prove the result for $p=1$. The case $p>1$ follows from the monotonicity of the average. Finally, we address the case $0<p<1$ through the same interpolation argument as in the proof of Proposition \ref{prop_lplinftycalpha}. 

In the sequel, we remove the assumption $u\in C^1(B_R)$ and complete the proof in the general setting.

\medskip

\noindent{\bf Step 4 - }Suppose $u\in C(B_R)$ satisfies \eqref{eq_belowc1a}, for some $P\in\mathbb{R}^d$. Next, we show that a mollification of $u$ also satisfies \eqref{eq_belowc1a}.

Fix $0<\varepsilon\ll1$. Let $x_0\in B_{R/4}$ and $y\in B_\varepsilon$, so that $x_0-y\in B_{R/4+\varepsilon}\subset B_{R/2}$. For $x\in B_{R/16}(x_0)$, we have $x-y\in B_{R/16}(x_0-y)$. Hence,
\begin{equation}\label{eq_regularize}
	u(x-y)\geq u(x_0-y)+\left\langle P,x-x_0\right\rangle -|x-x_0|\sigma(|x-x_0|).
\end{equation}

On the other hand, let $(\eta_\varepsilon)_{\varepsilon>0}$ be a family of standard symmetric mollifying kernels. Denote with $u_\varepsilon(x)$ the mollification of $u$ with respect to $\eta_\varepsilon$; that is,
\[
	u_\varepsilon(x):=\int_{B_\varepsilon}u(x-y)\eta_\varepsilon(y){\rm d}y.
\]
Multiplying \eqref{eq_regularize} by $\eta_\varepsilon$ and integrating over $B_\varepsilon$ we get
\[
	u_\varepsilon(x)\geq u_\varepsilon(x_0)+\left\langle P,x-x_0\right\rangle-|x-x_0|\sigma(|x-x_0|)
\]
for every $x_0\in B_{R/4}$ and $x\in B_{R/16}(x_0)$. 

Because $u_\varepsilon\in C^1(B_R)$, we have
\[
	\left\|u_\varepsilon\right\|_{L^\infty(B_{R/16})}\leq C\left(\intav{B_{R/16}}|u_\varepsilon(y)|{\rm d}y+R\sigma(R)\right),
\]
where the constant $C>0$ depends only on the dimension. By letting $\varepsilon\to 0$ one recovers the estimate for $u$ and completes the proof.
\end{proof}

We notice the proofs of Propositions \ref{prop_lplinftycalpha} and \ref{prop_lplinftyc1alpha} are inspired by ideas in \cite{Evans_Gariepy_1992}. In addition to their intrinsic interest, those results are useful in the analysis of several consequences of the Theorem \ref{thm_main1} and Corollary \ref{cor_main1} to the theory of nonlinear elliptic problems. In the sequel, we detail several such consequences.
	
\section{Consequences for nonlinear elliptic problems}\label{sec_cons}

In what follows we detail some consequences of the Theorem \ref{thm_main1} and Corollary \ref{cor_main1}. We start with the inclusion of functions satisfying {\rm [WH]} and [$L^p L^\infty$] in a class of viscosity solutions.

\subsection{Weak Harnack, $L^p L^\infty$-estimates and viscosity classes}\label{subsec_whleliimplyvisc}

Suppose $u\in C(B_1)$ belongs to the class $\overline{S}(\lambda,\Lambda,f)$, with $f\in L^\infty(B_1)\cap C(B_1)$. It is well-known that $u$ satisfies the weak Harnack inequality as well as an $L^p L^\infty$-estimate. Our next result establishes the opposite implication. That is, if a function satisfies {\rm [WH]} and [$L^p L^\infty$], then it is a supersolution to a fully nonlinear uniformly elliptic equation.

\begin{thm}\label{thm_whleliimplyvisc}
Let $u\in C(B_1)$. For any affine function $\ell(\cdot)$, suppose $u+\ell$ satisfies {\rm [WH]}, with $f\in L^\infty(B_1)$ and $\vartheta\equiv 2$. Suppose also that $u+\ell$ satisfies {\rm [$L^p L^\infty$]} with $\sigma(t):=At^2$, for some $A\geq 0$. Then
\[
	u\in\overline{\mathcal{S}}\left(L,L,A+\left\|f\right\|_{L^\infty(B_1)}\right),
\]
where $L=L(C_{\rm WH},C_{p,\infty})$.
\end{thm}
\begin{proof}
For the sake of presentation, we split the proof into three steps.

\medskip

\noindent{\bf Step 1 - }Suppose that $P(x)$ is a paraboloid touching the graph of $u$ from below at $x_0$ in $B_\delta(x_0)\subset B_1$. For simplicity, we suppose $x_0\equiv 0$ and write
\[
	P(x)=P(0)+DP(0)\cdot x+\left\langle D^2P(0)x,x\right\rangle.
\]
Set $\overline{A}:=\|(D^2P(0))^-\|$; we claim that 
\begin{equation}\label{eq_claim1}
	\left\langle D^2P(0)x,x\right\rangle\geq -\overline{A}|x|^2
\end{equation}
for every $x\in\mathbb{R}^d$. Indeed, by choosing an appropriate basis, we can write
\[
	\left\langle D^2P(0)x,x\right\rangle=\sum_{e_i>0}e_ix_i^2+\sum_{e_i<0}e_ix_i^2,
\]
for some $d$-uple $(e_1,\ldots,e_d)$. Thus
\begin{align*}
	\left\langle D^2P(0)x,x\right\rangle\geq-\sum_{e_i<0}|e_i|x_i^2\geq -\overline{A}\sum_{e_i<0}x_i^2\geq -\overline{A}|x|^2,
\end{align*}
which yields \eqref{eq_claim1}, verifying the claim.

\medskip

\noindent{\bf Step 2 - }By assumption, $u(x)\geq P(x)$ in $B_\delta$. Hence, 
\[
	u(x)-P(0)-DP(0)\cdot x\geq -\overline{A}|x|^2,
\]
for every $x\in B_\delta$. Denoting $\ell(x):=P(0)+DP(0)\cdot x$, we conclude that 
\[
	\inf_{B_\delta}\left(u-\ell\right)\geq -\overline{A}\delta^2.
\]
An application of Theorem \ref{thm_main1} yields
\begin{equation}\label{eq_oww}
	\sup_{B_{\delta/4}}\left(u-\ell\right)\leq C\left(A+\left\|f\right\|_{L^\infty(B_1)}+\overline{A}\right)\delta^2,
\end{equation}
where $C=C(C_{\rm WH},C_{p,\infty})$. Define $T(x)$ as
\[
	T(x):=C\left(A+\left\|f\right\|_{L^\infty(B_1)}+\overline{A}\right)|x|^2.
\]

\medskip

\noindent{\bf Step 3 - }From \eqref{eq_oww} we infer that $v:=u-\ell$ touches $T(x)$ from below at $x_0\equiv0$ in $B_{\delta/4}$. In addition, because $P$ touches $u$ from above at $x_0\equiv0$, we have
\[
	v(x)\geq \left\langle D^2P(0)x,x\right\rangle
\]
in $B_{\delta}$, with equality at $x=0$. We conclude that
\[
	\left\langle D^2P(0)x,x\right\rangle\geq C\left(A+\left\|f\right\|_{L^\infty(B_1)}+\overline{A}\right)|x|^2
\]
in $B_{\delta/4}$, with equality at $x=0$. The definition of $\overline{A}$ builds upon Lemma \ref{lem_curvature} to yield
\[
	\left\|\left(D^2P(x_0)\right)^+\right\|\leq C\left(A+\left\|f\right\|_{L^\infty(B_1)}+\left\|\left(D^2P(x_0)\right)^-\right\|\right)
\]
and complete the proof.
\end{proof}

\begin{Remark}[Sufficient conditions to become a solution]\label{rmk_solvability}\normalfont
In Theorem \ref{thm_whleliimplyvisc} we suppose that $u+\ell$ satisfies {\rm [WH]} and [$L^p L^\infty$] to conclude $u\in\overline{\mathcal{S}}$. By supposing that $-u+\ell$ satisfies {\rm [WH]} and [$L^p L^\infty$]  for any affine function $\ell$, one concludes $u\in\underline{\mathcal{S}}$. Hence, the proof of Theorem \ref{thm_whleliimplyvisc} yields the following: suppose $u\in C(B_1)$  is such that $(\pm u+\ell)$ satisfy {\rm [WH]} and [$L^p L^\infty$]; then $u$ is in the class $\mathcal{S}:=\overline{\mathcal{S}}\cap\underline{\mathcal{S}}$. C.f. \cite[Theorem 4]{Caffarelli_1999}.
\end{Remark}

\subsection{H\"older continuity for $\omega$-semiconvex supersolutions}\label{subsec_escauriazar}

In what follows we study $\omega$-semiconvex supersolutions to fully nonlinear equations. The distinctive feature concerns the integrability of the source term, which is strictly below the dimension and above the so-called Escauriaza's exponent.  

Let $u\in C(B_1)$ be an $\omega$-semiconvex viscosity solution to 
\begin{equation}\label{eq_dakota}
	\mathcal{M}^-(D^2u)\leq f\hspace{.2in}\mbox{in}\hspace{.2in}B_1,
\end{equation}
where $f\in L^q(B_1)$, for $p_0\leq q<d$. As a consequence of Corollary \ref{cor_main1}, we conclude that $u$ is locally H\"older-continuous, with estimates depending on the data of the problem.

\begin{thm}[H\"older continuity in the Escauriaza range]\label{thm_hcer}
Let $u\in \overline{S}(\lambda,\Lambda,f)$ in $B_1$. Suppose $u$ is an $\omega$-semiconvex function and $f\in L^q(B_1)$, for $p_0\leq q<d$. Then $u\in C^{2-d/q}_{\rm loc}(B_1)$. In addition, there exists a positive constant $C=C(d,q,\lambda,\Lambda,\omega)$ such that 
\[
	\left\|u\right\|_{C^{2-d/q}(B_{1/2})}\leq C\left(1+\omega(1)+\left\|u\right\|_{L^\infty(B_1)}+\left\|f\right\|_{L^q(B_1)}\right).
\]
\end{thm}
\begin{proof}
We start by noticing that translations of $u$ satisfy {\rm [WH]} and [$L^p L^\infty$]. Indeed, for any $c\in\mathbb{R}$, define $v:=u+c$. Clearly, $v\in \overline{S}(\lambda,\Lambda,f)$. Moreover, Remark \ref{rem_semiconvex} ensures that $v$ is $\omega$-semiconvex. Hence, because $v\in \overline{S}(\lambda,\Lambda,f)$, we infer it satisfies {\rm [WH]}. In addition, the $\omega$-semiconvexity of $v$ builds upon Lemma \ref{lem_BFMprop82b} to ensure that $v$ satisfies [$L^p L^\infty$]. In the sequel, we choose $c\equiv u(x_0)$ for an arbitrary $x_0\in B_{1/2}$, estimate some ingredients and resort to Corollary \ref{cor_main1}. For the sake of clarity, we split the proof into three steps.

\medskip

\noindent{\bf Step 1 - }The $\omega$-semiconvexity of $u$ ensures that 
\[
	u(x)-u(x_0)\geq \left\langle P,x-x_0\right\rangle-|x-x_0|\omega(|x-x_0|)
\]
for every $x\in B_1$, $P\in\partial_\omega u(x_0)$, and $x_0\in B_{1/2}$ fixed, though arbitrary. We continue by producing a modulus of continuity for $u(x)-u(x_0)$ from below. 

Let $\delta_0\in(0,1/4)$ be such that $\omega(\delta_0)\leq 1/2$. Define $\delta_1>0$ as
\[
	\delta_1:=\sup\{t\in(0,1/4)\;|\;\omega(t)\leq 1/2\};
\]
clearly, $0<\delta_0\leq\delta_1\leq 1/4$. For every $\rho\in(0,\delta_1)$ and $x\in B_\rho(x_0)$, we get
\begin{align}\label{eq_P1}\nonumber
	|\left\langle P,x-x_0\right\rangle-|x-x_0|\omega(|x-x_0|)|&\leq \big(|P|+\omega(|x-x_0|)\big)|x-x_0|\\\nonumber
		&\leq \big(|P|+\omega(\delta_1)\big)\rho\\
		&\leq \big(|P|+1\big)\rho.
\end{align}

On the other hand, because $P\in\partial_\omega u(x_0)$, Lemma \ref{lem_BMProp51} ensures that 
\begin{equation}\label{eq_P2}
	|P|\leq 8\left(\left\|u\right\|_{L^\infty(B_1)}+\omega(1)\right).
\end{equation}
Combining \eqref{eq_P1} and \eqref{eq_P2} one gets
\[
	|\left\langle P,x-x_0\right\rangle-|x-x_0|\omega(|x-x_0|)|\leq 8\left(\left\|u\right\|_{L^\infty(B_1)}+\omega(1)+1\right)\rho,
\]
and concludes 
\begin{equation}\label{eq_adriatico}
	\inf_{B_\rho(x_0)}\big(u-u(x_0)\big)\geq -8\left(\left\|u\right\|_{L^\infty(B_1)}+\omega(1)+1\right)\rho.
\end{equation}

\noindent{\bf Step 2 - }The discussion before Step 1 ensures that $v:=u-u(x_0)$ falls within the scope of the Corollary \ref{cor_main1}. As a consequence, the modulus of continuity in \eqref{eq_adriatico} yields 
\begin{equation}\label{eq_applyus}
	\sup_{B_{\rho/4}(x_0)}\big(u-u(x_0)\big)\leq C_1\left(1+\omega(1)+\left\|u\right\|_{L^\infty(B_1)}+\|f\|_{L^q(B_1)}\right)\rho^{2-\frac{d}{q}},
\end{equation}
for every $\rho\in(0,\delta_1)$, where $C_1=C_1(C_{\rm WH},C_{p,\infty})$. For simplicity, we define the constant $M$ as
\[
	M:=C_1\left(1+\omega(1)+\left\|u\right\|_{L^\infty(B_1)}+\|f\|_{L^q(B_1)}\right);
\]
hence \eqref{eq_applyus} becomes 
\[
	\sup_{B_{\rho/4}(x_0)}\big(u-u(x_0)\big)\leq M\rho^{2-\frac{d}{q}}.
\]

\medskip

\noindent{\bf Step 3 - }To complete the proof, let $x,y\in B_{1/2}$. Two cases arise; suppose first that $|x-y|\leq\delta_1$. By setting $\rho:=|x-y|$ we conclude
\[
	|u(x)-u(y)|\leq M\rho^{2-\frac{d}{q}}=M|x-y|^{2-\frac{d}{q}}.
\]
Conversely, if $|x-y|>\delta_1$, we obtain
\[
	\frac{|u(x)-u(y)|}{|x-y|^{2-\frac{d}{q}}}\leq \frac{2\|u\|_{L^\infty(B_1)}}{\delta_1^{2-\frac{d}{q}}}.
\]
In any case, one gets
\[
	|u(x)-u(y)|\leq\left(M+\frac{2\|u\|_{L^\infty(B_1)}}{\delta_1^{2-\frac{d}{q}}}\right)|x-y|^{2-\frac{d}{q}}
\]
and completes the proof.
\end{proof}

\subsection{Flipping H\"older continuity for supersolutions}\label{subsec_superholder}

In what follows we consider supersolutions $u\in \overline{S}(\lambda,\Lambda,f)$, with $f\in L^q(B_1)$ both for $p_0\leq q\leq d$ and $q>d$. However, instead of working under a (global) $\omega$-semiconvexity condition, we suppose $u$ has a (local) H\"older modulus of continuity \emph{from below}.

Therefore, we focus on supersolutions satisfying a (local) H\"older estimate from below. Our conclusion is that $u\in C^\alpha(B_1)$, with estimates. That is, the one-sided H\"older modulus of continuity becomes a two-sided regularity result, provided $u\in\overline{S}(\lambda,\Lambda,f)$. Our findings are the subject of the next theorem.

\begin{thm}\label{thm_holderfrombelow}
Fix $0<\lambda\leq\Lambda$. Let $u\in C(B_1)$ be a viscosity solution to
\[
	\mathcal{M}^-_{\lambda,\Lambda}(D^2u)\leq f\hspace{.3in}\mbox{in}\hspace{.3in}B_1,
\]
where $f\in L^q(B_1)$, for $p_0\leq q<d$. Suppose there exist $\alpha\in(0,1)$ and $C_\alpha>0$ such that
\begin{equation}\label{eq_sucessocortona}
	\inf_{B_\rho(x_0)}\big(u-u(x_0)\big)\geq - C_\alpha\rho^\alpha
\end{equation}
for every $x_0\in B_{1/2}$. Then $u\in C^\gamma_{\rm loc}(B_1)$, with 
\[
	\gamma:=\min\left\lbrace\alpha,2-\frac{d}{q} \right\rbrace,
\]
and there exists $C=C(d,q,\lambda,\Lambda,\alpha)$ such that
\[
	\|u\|_{C^\gamma(B_{1/2})}\leq C\left(C_\alpha+\|u\|_{L^\infty(B_1)}+\|f\|_{L^q(B_1)}\right).
\]
\end{thm}

A striking distinction between Theorem \ref{thm_holderfrombelow} and Theorem \ref{thm_hcer} sits in a local-global dichotomy. In fact, under an $\omega$-semiconvexity condition, the modulus of continuity satisfied by the supersolution is \emph{globally below} $u$. This is not the case under the assumptions of the Theorem \ref{thm_holderfrombelow}.  As a consequence, the $L^p L^\infty$-estimates recalled in Lemma \ref{lem_BFMprop82b} are not available, and we resort to Proposition \ref{prop_lplinftycalpha}. In what follows, we present the proof of the Theorem \ref{thm_holderfrombelow}.

\begin{proof}[Proof of Theorem \ref{thm_holderfrombelow}]
We aim at applying Corollary \ref{cor_main1} to $u-u(x_0)$, for any $x_0\in B_{1/2}$. Because 
\[
	\mathcal{M}^-_{\lambda,\Lambda}(D^2u)\leq f\hspace{.3in}\mbox{in}\hspace{.3in}B_1,
\]
we also have $v:=u-u(x_0)\in\overline{S}(\lambda,\Lambda,f)$. Hence, $v$ satisfies {\rm [WH]}. In addition,
\[
	\inf_{B_\rho(x_0)}\big(v-v(x_0)\big)=\inf_{B_\rho(x_0)}\big(u-u(x_0)\big)\geq - C_\alpha\rho^\alpha,
\]
and Proposition \ref{prop_lplinftycalpha} ensures that $v$ satisfies [$L^p L^\infty$]. An application of Corollary \ref{cor_main1} yields
\begin{equation}\label{eq_rome}
	\sup_{B_{\rho/4}(x_0)}\big(u-u(x_0)\big)\leq C\left(C_\alpha\rho^\alpha+\rho^{2-\frac{d}{q}}\|f\|_{L^q(B_1)}\right),
\end{equation}
for every $x_0\in B_{1/2}$, and every $\rho\in(0,1/8)$.  By combining \eqref{eq_sucessocortona} and \eqref{eq_rome}, we obtain
\begin{equation}\label{eq_conclusionh}
	\big|u(x)-u(x_0)\big|\leq C\left(C_\alpha+\|f\|_{L^q(B_1)}\right)|x-x_0|^\gamma
\end{equation}
provided $x_0\in B_{1/2}$, $x\in B_\rho(x_0)$, and $\rho\in (0,1/4)$. 

To complete the proof, we consider arbitrary $x,y\in B_{1/2}$ and examine two cases. First, suppose $|x-y|\leq 1/8$; hence, the estimate in \eqref{eq_conclusionh} is available. Suppose otherwise; that is, $|x-y|>1/8$. Then we have
\[
	\frac{\big|u(x)-u(y)\big|}{|x-y|^\gamma}\leq 2^{3\gamma+1}\|u\|_{L^\infty(B_1)}.
\]
In any case, we conclude
\[
	\big|u(x)-u(x_0)\big|\leq \big(2^{3\gamma+1}+C\big)\left(C_\alpha+\|u\|_{L^\infty(B_1)}+\|f\|_{L^q(B_1)}\right)|x-x_0|^\gamma
\]
and complete the proof.
\end{proof}

The result in Theorem \ref{thm_holderfrombelow} has a counterpart at the level of the gradient. I.e., if a supersolution $u:B_1\to\mathbb{R}$ has a $C^{1,\alpha}$-modulus of continuity from below, it is indeed $C^{1,\gamma}$-continuous, where the smoothness degree $\gamma\in(0,1)$ depends on the data of the problem. This is the content of the next theorem.

\begin{thm}\label{thm_flipc1a}
Fix $0<\lambda\leq\Lambda$. Let $u\in C(B_1)$ be a viscosity solution to
\[
	\mathcal{M}^-_{\lambda,\Lambda}(D^2u)\leq f\hspace{.3in}\mbox{in}\hspace{.3in}B_1,
\]
where $f\in L^q(B_1)$, for $q>d$. For every $x_0\in B_{1/2}$, suppose there exist an affine function $\ell_{x_0}(\cdot)$, $\alpha\in(0,1)$, and $C_\alpha>0$ such that
\begin{equation}\label{eq_sucessocortonab}
	\inf_{B_\rho(x_0)}\big(u-\ell_{x_0}\big)\geq - C_\alpha\rho^{1+\alpha},
\end{equation}
for every $\rho\in (0,1/2)$. Then $u\in C^{1+\gamma}_{\rm loc}(B_1)$, 
with 
\[
	\gamma:=\min\left\lbrace\alpha,1-\frac{d}{q} \right\rbrace.
\]
Moreover, there exists $C=C(d,q,\lambda,\Lambda,\alpha)$ such that
\[
	\|u\|_{C^{1+\gamma}(B_{1/64})}\leq C\left(C_\alpha+\|f\|_{L^q(B_1)}\right).
\]
\end{thm}
\begin{proof}
Because the extremal operator $\mathcal{M}^-_{\lambda,\Lambda}$ is invariant by affine translations, we have
\[
	\mathcal{M}^-_{\lambda,\Lambda}\big(D^2(u+\ell)\big)\leq f\hspace{.3in}\mbox{in}\hspace{.3in}B_1,	
\]
for every affine function $\ell$. Hence, $u+\ell$ satisfies {\rm [WH]} with $f$. In addition, Proposition \ref{prop_lplinftyc1alpha} ensures $u+\ell$ satisfies [$L^p L^\infty$]. Theorem \ref{thm_main1} implies
\begin{equation}\label{eq_sucessocortonanew}
	\sup_{B_{\rho/4}(x_0)}\big(u-\ell_{x_0}\big)\leq C\left(C_\alpha\rho^{1+\alpha}+\|f\|_{L^q(B_1)}\rho^{2-\frac{d}{q}}\right),
\end{equation}
where $C=C(d,p)$ is the same as in Proposition \ref{prop_lplinftyc1alpha}.

By combining \eqref{eq_sucessocortonab} with \eqref{eq_sucessocortonanew} one concludes there exists an affine function $\ell$, satisfying $\ell(x_0)=u(x_0)$, such that 
\begin{equation*}\label{eq_bologna}
	-C_1\rho^{1+\gamma}\leq \inf_{B_{\rho/4}(x_0)}\big(u(x)-\ell(x)\big)\leq \sup_{B_{\rho/4}(x_0)}\big(u(x)-\ell(x)\big)\leq C_1\rho^{1+\gamma},
\end{equation*}
for
\[
	C_1:=C\left(C_\alpha\rho^{1+\alpha}+\|f\|_{L^q(B_1)}\rho^{2-\frac{d}{q}}\right),
\]
where $C>0$ is universal. A straightforward application of Lemma 9.1 in \cite{Braga_Figalli_Moreira_2019} completes the proof.
\end{proof}

We close this section with a discussion on the so-called De Giorgi class. It comprises a wide latitude of functions: from solutions to elliptic equations in divergence form to minimizers and $Q$-minimizers of variational integrals.

\subsection{De Giorgi classes}\label{subsec_degiorgi}

In what follows we mention the De Giorgi class as a collection of functions satisfying {\rm [WH]} and [$L^p L^\infty$]. The De Giorgi class is defined in terms of energy inequalities on super-level sets. To be more precise, we proceed with a definition.

\begin{Definition}[De Giorgi classes]\label{def_degiorgiclass}
We say that $u\in W^{1,p}_{\rm loc}(\Omega)$ is in the De Giorgi class. be an open set. A function $u\in W^{1,p}_{\rm loc}(\Omega)$ is in the De Giorgi class ${\rm DG}_p^\pm(\Omega)$ if, for every $B_R(y)\subset\Omega$, $0<\sigma<1$, and $k>0$, we have
\begin{align}\label{eq_degiorgicac}
	\int_{B_{\sigma R}(y)}|D(u-k)^\pm|^p{\rm d}x&\leq \frac{C_1}{(1-\sigma)^pR^p}\int_{B_R(y)}|(u-k)^\pm|^p{\rm d}x\\\nonumber
		&\quad+C_1\left[\chi^p+\left(R^{-\frac{d\varepsilon}{p}}k\right)^p\right]|A_{k,R}^\pm|^{1-\frac{p}{d}+\varepsilon},
\end{align}
where $C_1>0$, $\chi>0$, and $\varepsilon\in (0,p/d)$ are constants, and 
\[
	A_{k,R}^\pm:=\left\lbrace x\in B_R(y)\,|\,(u-k)^\pm>0\right\rbrace.
\]
We define ${\rm DG}_p(\Omega):={\rm DG}_p^+(\Omega)\cap{\rm DG}_p^-(\Omega)$.
\end{Definition}

In \cite{DiBenedetto_Trudinger_1984} the authors prove that functions in the De Giorgi class ${\rm DG}_p^+(\Omega)$ satisfy {\rm [WH]} and [$L^p L^\infty$]. As a consequence, those functions fall within the scope of the Corollary \ref{cor_main1}. We conclude that functions in the ${\rm DG}_p^+(\Omega)$ are entitled to the conclusions of Theorems \ref{thm_hcer} and \ref{thm_holderfrombelow}. Rigorously put, one concludes:

\begin{Proposition}[De Giorgi class]\label{prop_dgc}
Let $\Omega\subset\mathbb{R}^d$ and suppose $(u+c)\in {\rm DG}_p^-(\Omega)$, for every $c\in\mathbb{R}$. Suppose further there exists a modulus of continuity $\gamma(\cdot)$ such that
\[
	\inf_{x\in B_\rho(x_0)} \big(u-u(x_0)\big)\geq -\gamma(\rho),
\]
for every $x_0\in B_{1/2}$ and $\rho\in(0,1/4)$. Then
\[
	\sup_{B_{\rho/4}(x_0)} \left(u-u(x_0)\right)\leq C\left(\chi\rho^\frac{d\varepsilon}{p}+\gamma(\rho)\right),
\]
where $C=C(d,p,C_1,\varepsilon)$.
\end{Proposition}
\begin{proof}
The result follows directly from Corollary \ref{cor_main1} combined with Proposition \ref{prop_lplinftycalpha} and the Theorem 2 in \cite{DiBenedetto_Trudinger_1984}.
\end{proof}

\bigskip

\noindent{\bf Acknowledgements.} Part of this work was developed during the authors' stay at The Abdus Salam International Centre for Theoretical Physics (ICTP); the authors thank ICTP for its hospitality and the vibrant scientific environment. DM is partly funded by CNPq and FUNCAP (PRONEX). EP is partly funded by FAPERJ (E-26/201.390/2021) and the Centre for Mathematics of the University of Coimbra. This work was partially supported by the Centre for Mathematics of the University of Coimbra - UIDB/00324/2020, funded by the Portuguese Government through FCT/MCTES.

\bibliography{biblio}
\bibliographystyle{plain}

\bigskip
\noindent\textsc{Diego R. Moreira}\\Departamento de Matem\'atica\\Universidade Federal do Cear\'a (Fortaleza, Brazil)\\\noindent\texttt{dmoreira@mat.ufc.br}

\bigskip

\noindent\textsc{Edgard A. Pimentel}\\University of Coimbra\\CMUC, Department of Mathematics\\ 3001-501 Coimbra, Portugal\\and\\Pontifical Catholic University of Rio de Janeiro -- PUC-Rio\\22451-900, G\'avea, Rio de Janeiro-RJ, Brazil\\\noindent\texttt{edgard.pimentel@mat.uc.pt}

\end{document}